\documentclass[11pt,a4paper]{article}
\usepackage[]{geometry}
\usepackage[mathcal]{euscript}
\usepackage{bm,url}                     
\usepackage{amsfonts}
\usepackage{amssymb}
\usepackage{amsmath}
\usepackage{amsthm}
\usepackage[utf8x]{inputenc}
\usepackage{float}
\usepackage[pdftex]{graphicx}
\DeclareGraphicsExtensions{.bmp,.png,.pdf,.jpg,.ps}
\usepackage{colortbl}  
\usepackage{todonotes}             
\usepackage{booktabs}
\newtheorem{definition}{Definition}
\newtheorem{remark}{Remark}
\newtheorem{theorem}{Theorem}

\newtheorem{notation}[theorem]{Notation}

\usepackage[english]{babel}
\usepackage{amsmath}
\usepackage{amsfonts}
\usepackage{amssymb}
\usepackage{graphicx}
\usepackage{color}
\usepackage{array}
\usepackage{multirow}
\usepackage{mathrsfs}
\usepackage{hyperref}

\newcommand{\pr}{\mathbb{P}}

\definecolor{violet}{rgb}{0.7,0,0.6}
\definecolor{OliveGreen}{RGB}{85,107,47}

\usepackage{xcolor}
\newif\ifnotes\notestrue
\title{Short range vs long range dependence. An hypothesis test based on  Fractional Iterated
	Ornstein--Uhlenbeck processes}
\author{Juan Kalemkerian \\
Universidad de la República, Facultad de Ciencias\\ \\
Andr\'es Sosa\\
Universidad de la República, Facultad de Ciencias Económicas \\ y Administración}
\begin{document}
\maketitle

\begin{abstract}
\noindent  In this work, which is based on the family of Fractional Iterated Ornstein--Uhlenbeck processes, we propose a new hypothesis test to contrast short memory versus long memory in time series. This family includes short memory and
long memory processes, and has the ability to approximate a long memory processes by a short memory processes. Based on the asymptotic results of the estimators of its parameters, we will  present the test and show how it can be implemented. Also, we will show a comparison with other tests widely used under both short memory 
and long memory scenarios. The main conclusion is that this new test is the one with best performance under the null hypothesis, and has the maximum power in some alternatives.

\end{abstract}


  \section{Introduction} 
  It is commonly said that the sequence has long range dependence when we have a stationary and centered sequence of random variables $X_1,X_2,....,X_n,...$ such that 
  $\mathbb{E}(X_0X_n) \rightarrow 0$ and $\sum_{n=1}^{+\infty}\left | \mathbb{E}(X_0X_n)\right |=+\infty$. In contrast, the sequence has short range dependence
  when it is satisfied that 
  $\sum_{n=1}^{+\infty}\left | \mathbb{E}(X_0X_n)\right |<+\infty$.
  Long memory was discovered by Hurst in his pioneering work (\cite{Hurst}), in which it shows its presence in a time series of variables (from several disciplines), such  as rainfall, temperature, pressure,  
  thickness of the rings of certain trees, sunspots and stock market phenomena. From that moment, long memory processes have been (and are still) studied extensively from both 
  theoretical and practical points of view.
  Although there are alternative definitions of short range or long
  range dependence (see for instance\cite{Brockwell}, \cite{Haslett}, \cite{Hipel}, \cite{Palma}), they are not all equivalent to each other. However, the underlying idea in all cases is that when the autocorrelation 
  function tends to zero quickly (slowly) we have a short (long) range process.  It is worth noting that there are processes to model both short memory and long memory time series, even though 
  it is clear that the long memory processes (such as ARFIMA process) are more complex models, are less intuitive and are more difficult to estimate.
  
  In empirical applications, it is important to have an hypothesis
  test that can help us to decide if we are dealing with a short memory or long memory time series. Unfortunately, there are not many hypothesis tests to tackle these kinds 
  of problems, while the existing tests have problems in their applicability. For example, they have a parameter in such way that the test can have a strong bias to reject the null hypothesis
  for some values of the parameter or a strong bias to non-reject the null hypothesis in the other cases.  This fact is not surprising due to the complexity of the problem to be solved. 
  Furthermore, even in those cases in which we take a value of the parameter where the test has a bias to  non-reject the null hypothesis, we will show in Section 5 that there are several 
  examples of short memory processes that are wrongly considered  (with very high probability) as long memory processes. Consequently, we can conclude that the existing 
  tests are not very reliable when working with real time series. This is an important problem because if we apply any of these hypothesis tests when we have a real time series, we cannot have 
  much confidence that we are making the correct decision. From a simple generalisation of the R/S statistic proposed by Hurst (\cite{Hurst}) in 1991, Lo (\cite{Lo}) develops a test that takes
  as the null hypothesis that the series has short memory dependence against the alternative that has long memory dependence. Although this is not the first hypothesis  
  test that has been developed for this purpose (see for example \cite{Geweke}), it is a widely used test. From this work, different variants have been developed that have generated other hypothesis tests, such as 
  Giraitis et al.'s (\cite{giraitis}) test that is based on the rescaled variance statistic, Lee and Schmidt's (\cite{Lee}) test, and others. 
  A description of these and other tests can be found in \cite{Beran}.

  Recently,
  a fractional iterated Ornstein--Uhlenbeck process was introduced in \cite{chichi}. 
  When the number of iterations in this family of processes is greater than or equal to $2$, they are short memory processes. However, in some cases we
  can approximate (in a continuous way) the fractional Ornstein--Uhlenbeck long range process. In addition, if we use only one iteration, then we have a fractional Ornstein--Uhlenbeck process. 
  This property can be use to design a hypothesis test to deal with this issue.  In this work, we propose a new hypothesis test that is based on a fractional iterated Ornstein--Uhlenbeck process
  to test short range dependence against long range dependence in observed  time series. In addition, we will make a comparison with other tests that are commonly used in practice using a wide 
  spectrum of scenarios, covering both short memory and long memory time series.

The rest of this paper is organised as follows.  In Section \ref{section2},  we  introduce the Fractional Iterated Ornstein--Uhlenbeck process, 
we explain the main properties that these processes satisfy  and we also show the main tools for the approach of the proposed  test. In Section \ref{section3}, 
we provide the test approach, which is based on the properties given in Section \ref{section2}. In Section \ref{section4}, we explain the hypothesis test implementation. In Section 
  \ref{section5}, we  compare the performance of this new test against other existing tests in the literature. In this section, we also include a criterion to select the parameter before
  using each of the competitor tests.
  In Section \ref{section6}, we  show an application to real data. Our concluding remarks are given in Section \ref{section7}.

\section{Preliminaries} \label{section2}
We start by defining the  fractional iterated Ornstein--Uhlenbeck processes of order $2$. 
We also summarise the main properties that will allow us to develop the idea of the hypothesis test.
An introduction to this type of process, as well as its theoretical development, can be found in 
  \cite{Kalemkerian_discreto} and \cite{chichi}. 

\begin{definition} \label{definicion1}
	Given $\ $a fractional Brownian motion $\left\{ B_{H}(t)\right\} _{t\in 
		\mathbb{R}}$,  a fractional Ornstein--Uhlenbeck process with
	parameters $\sigma ,\lambda >0$ and $H\in \left(  0,1\right] $  is defined as $\left\{
	X_{t}\right\} _{t\in \mathbb{R}}$ where 
	\begin{equation*}
	X_{t}=\sigma \int_{-\infty }^{t}e^{-\lambda \left( t-s\right) }dB_{H}(s)%
	\text{ for all }t\in \mathbb{R}.
	\end{equation*}
\end{definition}

In  \cite{Cheridito}, it is proved that the process  in Definition \ref{definicion1}  is the only stationary solution of 
the stochastic equation $$dX_t=-\lambda X_tdt+\sigma dB_H(t).$$ 
\begin{notation}
In this work, we use the  notation
	$\left\{ X_{t}\right\} _{t\in \mathbb{R}}$ $\sim $FOU$\left( \lambda ,\sigma
	,H\right) $ or  FOU $\left( \lambda,H\right) $ when $\sigma=1$.
\end{notation}

\begin{definition}
	Given $\ 0\leq \lambda _{1}<\lambda _{2},$ a fractional Brownian motion $%
	\left\{ B_{H}(t)\right\} _{t\in \mathbb{R}}$ and the processes  $\left\{ X_{t}^{\left(
		i\right) }\right\} _{t\in \mathbb{R}}$ that satisfying $X_{t}^{\left( i\right)
	}=\sigma \int_{-\infty }^{t}e^{-\lambda _{i}\left( t-s\right) }dB_{H}(s)$
	for $i=1,2,$. The fractional  iterated Ornstein--Uhlenbeck process $\left\{
	X_{t}\right\} _{t\in \mathbb{R}}$  with
	parameters $\lambda_1,\lambda_2$, $\sigma >0$ and $H\in \left( 0,1\right] $  is  
	\begin{equation*}
	X_{t}=\sigma \int_{-\infty }^{t}e^{-\lambda _{2}\left( t-s\right)
	}dX_{s}^{\left( 1\right) }\text{ for all }t\in \mathbb{R}.
	\end{equation*}
\end{definition}

\begin{notation}
Again, in this work we use the notation 
	$\left\{ X_{t}\right\} _{t\in \mathbb{R}}$ $\sim $FOU$\left( \lambda
	_{1},\lambda _{2},\sigma ,H\right) $ or  FOU
	 $\left( \lambda_1,\lambda_2,H\right) $ when $\sigma=1$.
\end{notation}

\begin{remark}
It is important to note that in the case that $\lambda_1=0$, we have the classic Ornstein-Uhlenbeck process.
\end{remark}

\begin{remark}
Observe that if we define the family of operators
 $$T_{\lambda}(y)(t)=\int_{-\infty}^{t}e^{-\lambda (t-s)}dy(s),$$
  therefore $X_t^{(i)}=T_{\lambda_i}(B_H)(t)$ for $i=1,2$, and 
$X_t=T_{\lambda_2}(X_.^{(1)})(t)=T_{\lambda_2} \left ( T_{\lambda_1}(B_H) ) \right (t)$. This properties enables us to state that any $FOU(\lambda_1,\lambda_2,H)$ is the composition of the operators $T_{\lambda_1}$ 
and $T_{\lambda_2}$ evaluated on a fractional Brownian motion with parameter $H$.
\end{remark}

\begin{remark}
	Integrating by parts, we can obtain that any $\left\{ X_{t}\right\} _{t\in 
		\mathbb{R}}$ $\sim $FOU$\left( \lambda _{1},\lambda _{2},\sigma ,H\right) $
	can be expressed by 
	\begin{equation*}
	X_{t}=\frac{\lambda _{1}}{\lambda _{1}-\lambda _{2}}X_{t}^{\left( 1\right) }+%
	\frac{\lambda _{2}}{\lambda _{2}-\lambda _{1}}X_{t}^{\left( 2\right) }
	\text{ for all }t\in \mathbb{R}.
	\end{equation*}
	Thus, any FOU$(\lambda_1,\lambda_2,H)$ is a linear combination of a FOU$(\lambda_1,H)$ and FOU$(\lambda_2,H)$ driven for the same 
	fractional Brownian motion and therefore the composition
	is commutative: $T_{\lambda_1} \left ( T_{\lambda_2}\right )=T_{\lambda_2} \left ( T_{\lambda_1} \right ).$
\end{remark}
\begin{remark}
We can generalise and compose $p$ times with the operators $T_\lambda$, by the same or different
values of $\lambda$ and obtain the FOU$(p)$ processes (according with notation given in \cite{chichi})
($p$ iterations). However, for our purposes it is enough to work with 
the FOU$(2)$ processes that is, FOU$(\lambda_1,\lambda_2, \sigma, H)$ where $0 \leq \lambda_1 < \lambda_2$.
\end{remark}

\begin{remark}
	Any $\left\{ X_{t}\right\} _{t\in \mathbb{R}}$ $\sim $FOU$\left( \lambda
	_{1},\lambda _{2},\sigma ,H\right) $ is a stationary centered Gaussian
	process.
\end{remark}

In the  memory process problem, it is important to know the autocorrelation function of the models. In the case of fractional Brownian motion $\left\{ X_{t}\right\} _{t\in \mathbb{R}}$ $%
\sim $FOU$\left( \lambda _{1},\lambda _{2},\sigma ,H\right) $  the
autocovariance function is given by 
\begin{equation}
\mathbb{E}\left( X_{0}X_{t}\right) =\frac{\sigma ^{2}H}{2}\left( \frac{\lambda
	_{1}^{2-2H}}{\lambda _{1}^{2}-\lambda _{2}^{2}}f_{H}\left( \lambda
_{1}t\right) -\frac{\lambda _{2}^{2-2H}}{\lambda _{1}^{2}-\lambda _{2}^{2}}%
f_{H}\left( \lambda _{2}t\right) \right), \label{cov_linear_combination}
\end{equation}%
where the function $f_{H}$ is defined by 
\begin{equation*}
f_{H}(t)=e^{-x}\left( \Gamma \left( 2H\right) -\int_{0}^{x
}e^{s}x^{s-1}dx\right) +e^{x}\left( \Gamma \left( 2H\right)
-\int_{0}^{x }e^{-s}x^{s-1}dx\right). 
\end{equation*}
For further details, the properties of the function $f_{H}$ can be found in \cite{chichi}. Between them, the asymptotic behavior of $f_H$: 
if  $H\neq 1/2$ we have that $f_{H}(x)\sim 2\left( 2H-1\right) x^{2H-2}$ as $x\rightarrow +\infty.$
In particular from (\ref{cov_linear_combination}), we have that the autocovariance function of any $%
\left\{ X_{t}\right\} _{t\in \mathbb{R}}$ $\sim $FOU$\left( \lambda ,\sigma
,H\right) $ is given by 
\begin{equation}
\mathbb{E}\left( X_{0}X_{t}\right) =\frac{\sigma ^{2}H}{2\lambda^{2H}}%
f_{H}\left( \lambda t\right). \label{covariance_fou}
\end{equation}
Therefore, the equation (\ref{covariance_fou})  allows us to conclude that  any FOU$%
\left( \lambda ,\sigma ,H\right) $ has short range dependence when $H\leq
1/2 $ and long range dependence for $H>1/2.$ 

Although from equation  (\ref{cov_linear_combination}) the autocovariance function of any FOU$\left( \lambda _{1},\lambda
_{2},\sigma ,H\right) $ is a linear combination of the functions $%
f_{H}\left( \lambda _{1}t\right) $ and $f_{H}\left( \lambda _{2}t\right).$ In addition,
 when $\lambda_1 >0$ any FOU$\left( \lambda _{1},\lambda _{2},\sigma
,H\right) $ process  has short range dependence for any value of $H\in \left(
0,1\right] .$ However, if $\lambda_1=0$, then we have that a FOU$(\lambda_2,\sigma,H)$  has short 
range dependence for $H \leq 1/2$ and long range dependence for $H>1/2$.

When we have a FOU$\left( \lambda _{1},\lambda _{2},\sigma ,H\right) $
process observed in an equispaced sample of the interval $\left[ 0,T\right]
, $ under some conditions that imply  $T\rightarrow +\infty $ and $T/n \rightarrow 0$, the parameters $\sigma ,H$ can be estimated
in a consistent way by a procedure proposed in the work \cite{Istas}. In Section \ref{section4} we
summarise the explicit
formulas
to obtain $\hat{H}$ and $\hat{\sigma}$ and their asymptotic behavior. 
These results are the basis for the  proposed hypothesis test. Moreover, it is satisfied that
 $\widehat{H}$ and $\widehat{\sigma }$  are asymptotically normal.

The structure of the spectral density is another interesting result  of fractional 
iterated Ornstein--Uhlenbeck processes. It is possible to prove 
that for each process FOU$\left( \lambda _{1},\lambda _{2},\sigma
,H\right),$ the spectral density    is given by 
\begin{equation}
f^{\left( X\right) }(x)=\frac{\sigma ^{2}\Gamma \left( 2H+1\right) \sin
	\left( H\pi \right) \left\vert x\right\vert ^{3-2H}}{2\pi \left( \lambda
	_{1}^{2}+x^{2}\right) \left( \lambda _{2}^{2}+x^{2}\right) }.
\label{spectral_density}
\end{equation}
Regarding the parameters $\lambda _{1}$ and $\lambda _{2}$, taking advantage
of the spectral density in the equation (\ref{spectral_density}), when $T\rightarrow +\infty $ can be
estimated in a consistent way by a modified Whittle contrast if the
process is observed in an equispaced sample of $\left[ 0,T\right]$ and has
asymptotically joint Gaussian distribution. 
In Section \ref{section4} we summarise the procedure to estimate
 $\hat{\lambda_1}$ and $\hat{\lambda_2}$. 
 \begin{remark}\label{H_interpretation}
  The interpretation of the $H$ parameter in any FOU$(\lambda_1,\lambda_2,\sigma,H)$ is related
  to the irregularity of the trajectories because $2H$ is the local  H\"older index of the process. 
  In the particular case in which
  $\lambda_1=0$, $H$ is also a parameter that govern the memory of the process (i.e., long  memory for $H > 1/2$ 
  and short memory
  for $H\leq 1/2$) in a similar way that it is interpreted in the fractional Gaussian noise or $d$ 
  parameter in ARFIMA model.
 \end{remark}

\section{Statistical hypothesis testing} \label{section3}
We assume that we have a sample $X_{t_{1}},X_{t_{2}},...,X_{t_{n}}$ of some centered stationary
process 
$\left\{ X_{t}\right\} _{t\in \mathbb{R}}$
where $%
0\leq t_{1}<t_{2}<...<t_{n}\leq T.$
Our  objective is to detect if there is a short range or long range dependence in the time series. Therefore,
we want to test 
\begin{eqnarray*}
	H_{0} &:&\left\{ X_{t}\right\} _{t\in \mathbb{R}}\text{has short range
		dependence} \\
	H_{1} &:&\left\{ X_{t}\right\} _{t\in \mathbb{R}}\text{has long range
		dependence.}
\end{eqnarray*}%

We  assume that the sample corresponds to
 $\left\{
X_{t}\right\} _{t\in \mathbb{R}}\sim $FOU$\left( \lambda _{1},\lambda
_{2}, \sigma, H\right) $ where $0\leq \lambda _{1}<\lambda _{2},$ $\sigma>0$, $H\in \left( 0,1%
\right] ,$ and the observations are equispaced. Therefore, according to  Section \ref{section2},  we can express the hypotheses test
in a parametric form as

\begin{eqnarray*}
	H_{0} &:&\text{ }H\leq 1/2\text{ \ or }\lambda _{1}>0 \\
	H_{1} &:&\text{ }H>1/2\text{ and }\lambda _{1}=0.
\end{eqnarray*}

Since that we have a consistent procedure to estimate $H$ and $\lambda _{1}$ (see Section \ref{section4}), it
is natural to reject the null hypothesis when $\widehat{H}\geq k$ and $%
\widehat{\lambda }_{1}\leq c$ where the values $k$ \ and $c$ are real constants such that 
\begin{equation*}
\sup_{H_{0}} \pr\left( \left\{ \widehat{H}\geq k\right\} \cap \left\{ \widehat{%
	\lambda }_{1}\leq c\right\} \right) =\alpha
\end{equation*}%
where $\alpha $ is the  signification level of the test. A simple way to obtain values of $k$ and $c$ such that the level of the test
is less than or equal to $\alpha, $ is to obtain $k$ such that $%
\sup_{H_{0}}\pr\left( \widehat{H}\geq k\right) =\alpha $ and $c$ such that $%
\sup_{H_{0}}\pr\left( \widehat{\lambda }_{1}\leq c\right) =\alpha .$

The asymptotic distribution of $\widehat{H}$ is normal and  is independent of the values of 
$\lambda _{1}$ and $\lambda _{2},$ then $$\sup_{H_{0}}\pr\left( \widehat{H}\geq
k\right) =\sup_{H\leq 1/2}\pr\left( \widehat{H}\geq k\right) $$
 and the supreme naturally
is reached at $H=1/2.$ Therefore, it is easy to find $k$ such that $%
\sup_{H_{0}}\pr\left( \widehat{H}\geq k\right) =\alpha .$ To obtain $\sup_{H_{0}}\pr\left( \widehat{\lambda }_{1}\leq c\right) $, we
observe that $\pr\left( \widehat{\lambda }_{1}\leq c\right) =g_{c}\left(
\lambda _{1},\lambda _{2},\sigma, H\right) $, and it is natural to expect that the
supreme is reached at the point (or points) where it is more difficult to decide whether 
$H_{0}$ or $H_{1}$ is true.

\begin{figure}[H]
 \centering
   \includegraphics[scale=0.5]{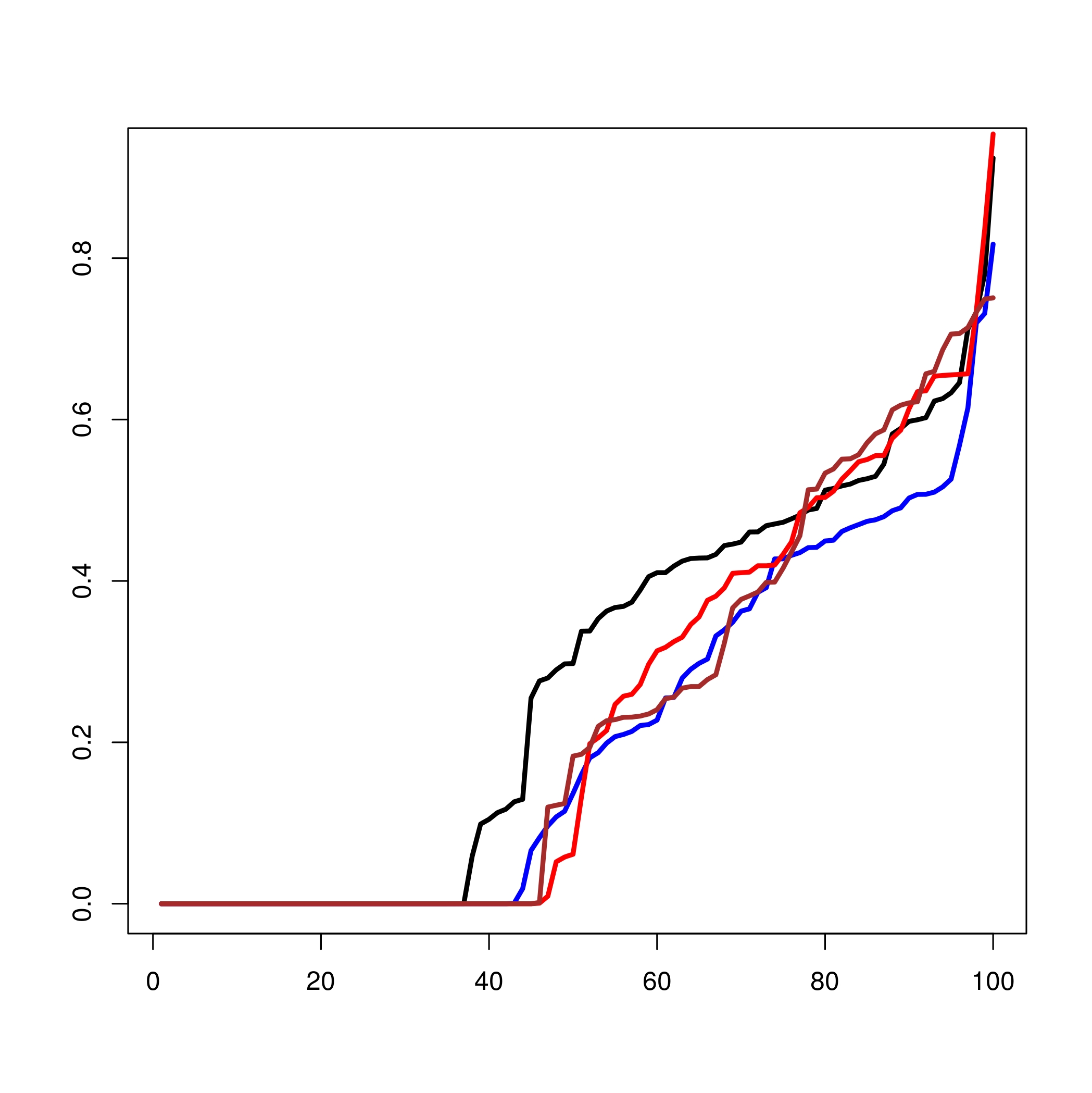} 
   \caption{$100$ values of $\hat{\lambda}_1$ for an equispaced sample of length $n=1000$ in $[0,100]$ of 
   a FOU$(\lambda_1=0.3,\lambda_2=0.8,\sigma=1, H)$ process for $H=0.2$ (black), $H=0.3$ (blue),
   $H=0.4$ (red) and $H=0.5$ (brown).
    } 
  \label{l1g}
\end{figure}

It is important to note that for small values of $c$, the function $%
g_{c}\left( \lambda _{1},\lambda _{2},\sigma, H\right) $ is increasing as an $H-$%
function for fixed values of $\sigma, \lambda_1$ and $\lambda_2$. For example,
for each $H=0.2, 0.3, 0.4, 0.5$ and $\lambda _{1}=0.3,$ $\lambda _{2}=0.8$,
$\sigma=1$,
we  simulate $100$ trajectories 
of FOU$(\lambda_1,\lambda_2, \sigma, H)$ 
in $[0,T]$ for $T=100$ and sample size of $n=1000$. In each case we have calculated $\hat{\lambda}_1$. 
In Figure \ref{l1g}, we show  $100$  values of 
$\hat{\lambda}_1$ ordered from the smallest to the largest.
For other values of $\lambda _{1},\lambda _{2}, \sigma$ the behavior is
similar.

Then, $\sup_{H_0}g_{c}\left( \lambda _{1},\lambda _{2},\sigma, H\right)
=g_{c}\left( \lambda _{1},\lambda _{2},\sigma, 1/2\right) .$

\begin{figure}[H]
 \centering
   \includegraphics[scale=0.5]{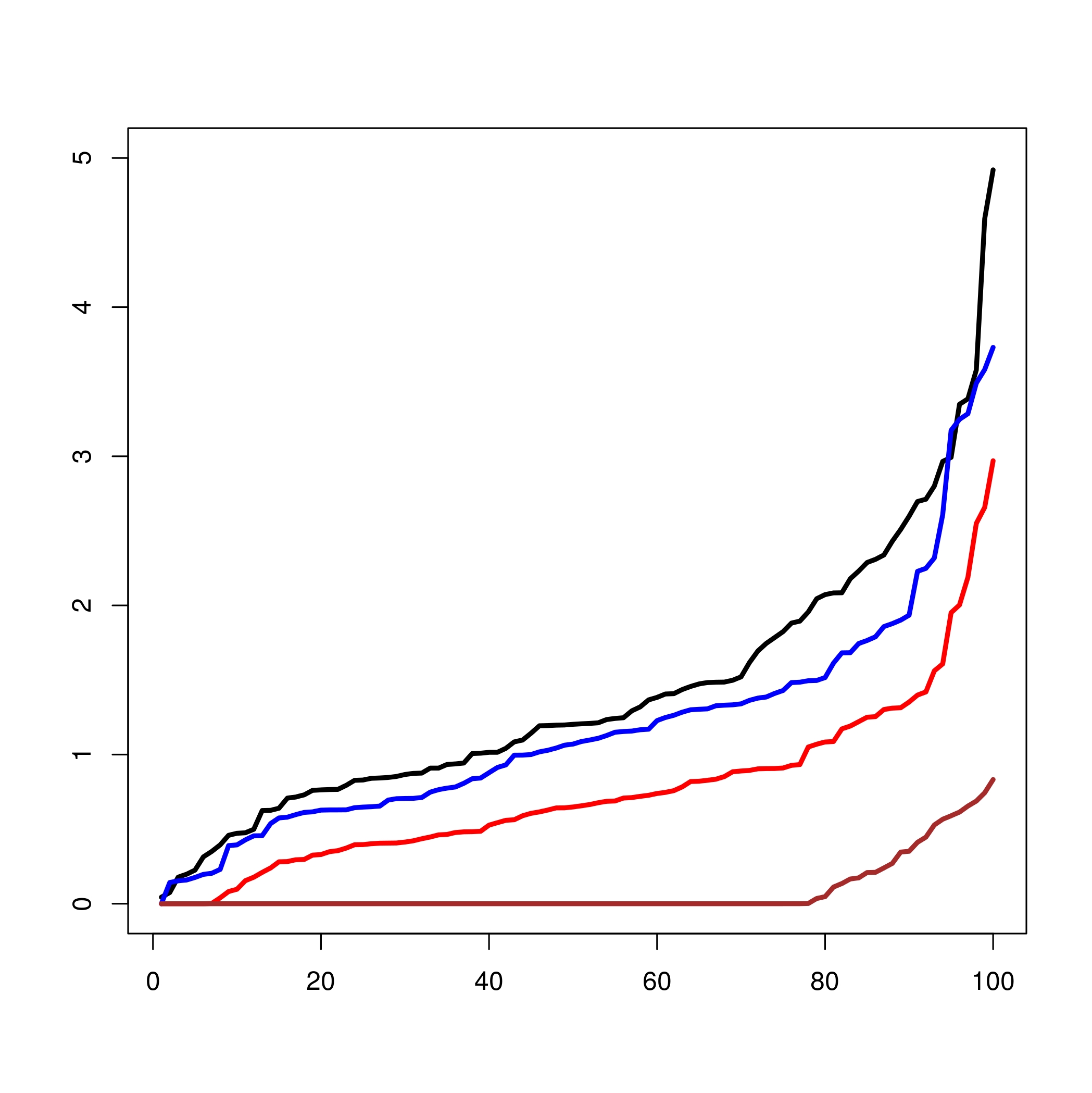} 
   \caption{$100$ values of $\hat{\lambda}_1$ for an equispaced sample of length $n=1000$ in $[0,100]$ of 
   a FOU$(\lambda_1,\lambda_2=5,\sigma=1, H=0.5)$ process for $\lambda_1=3$ (black), $\lambda_1=2$ (blue),
   $\lambda_1=1$ (red) and $\lambda_1=0$ (brown).
    } 
  \label{l1gH05}
\end{figure}

It should also be noted that the function $g_{c}\left( \lambda _{1},\lambda
_{2},\sigma, 1/2\right) $ grows as $\lambda _{2}-\lambda _{1}$ grows for moderately
small $c$ values. 
For example when $\lambda _{2}=5,$ $H=0.5$  and each $\lambda_1=0,1,2,3$ 
we have simulated $100$ 
process of FOU$(\lambda_1,\lambda_2, H)$ 
in $[0,T]$ for $T=100$ and sample size of $n=1000$.
In each case we have calculated $\hat{\lambda}_1$. 
In Figure \ref{l1gH05}, we show  the $100$  values of 
$\hat{\lambda}_1$ ordered from the smallest to the largest.
Again, the behavior is similar for other values of $\lambda _{1},\lambda _{2}$.

Therefore

\begin{equation*}
\sup_{H_{0}}P\left( \widehat{\lambda }_{1}\leq c\right) =\sup_{0<\lambda
	_{1}<\lambda _{2}\leq \widetilde{\lambda }}\sup_{H\leq 1/2}g_{c}\left(
\lambda _{1},\lambda _{2},\sigma,H\right) =g_{c}\left( 0,\widetilde{\lambda },%
\sigma,1/2\right) .
\end{equation*}
Taking into account that $\sigma $ does not appear in the hypotheses and that the behavior of $g_c$ as a function
of $H$ and as a function of $\lambda_1, \lambda_2$ is the same for any 
value of $\sigma$, we propose to obtain an approximated  value of $c$  using $\hat{\sigma}$ instead of $\sigma$, which is simply performed $$g_{c}\left( 0,\widetilde{\lambda },%
\widehat{\sigma},1/2\right)=\alpha.$$
\begin{remark}
 According with Remark \ref{H_interpretation}, we can have a long memory process where the H\"older index is
 less than $1/2$. For this case, the FOU 
 test naturally non rejects the null hypothesis. We observe examples of this situation in Section \ref{section5}. \label{H_reject}
\end{remark}

\section{Implementation of the test} \label{section4}
Given $X_1, X_2,...,X_n$ observations of a stationary centered time series, to perform the test we need to consider the observations as an equispaced sample of 
FOU$(\lambda_1,\lambda_2,\sigma, H)$ in some interval $[0,T]$; that is, $X_{T/n},X_{2T/n},...,X_{T}$.
To perform a hypothesis test, we need to estimate the parameters $\sigma$, $H$ and $\lambda_1$ which depend on
 $T$ (except $H$ as we will look in the following subsection). Thus,  we first need to know
the value of $T$. In Subsection 4.2 we propose a criterion to select a suitable value of $T$.

\subsection{Estimation of $\sigma, H$ and $\lambda_1$.}
If we know the value of $T$, we can proceed to estimate the parameters $\sigma, H$ and $\lambda_1$. The estimation of $H$ is independent of the value of $T$, the estimation of $\sigma$ depends on
$\widehat{H}$ and $T$, and the estimation of $\lambda_1$ requires knowledge the values of $ \widehat{\sigma}, \widehat{H}$ and $T$. More explicitly, the estimation of 
$\sigma, H$ and $\lambda_1$ is carried out in the  following two steps.\\

\noindent \textbf{Step 1. Estimation of $H$ and $\sigma$.}

First, we need to select a filter $a=(a_0,a_1,...,a_k)$ of length $k+1$ and order $L \geq 2$ (i.e.,
$\sum_{i=0}^{k}i^{j}a_i=0$ for $j < L$ and 
$\sum_{i=0}^{k}i^{L}a_i=0$).  For example we can  use $a=(-1,2,-1)$ (filter of order $2$) 
or the Daubechies filter of order $2$. given by 
 $$a=(0.4829629131445,-0.8365163037378,0.2241438680420,0.1294095225512)/\sqrt{2},$$
In this way,   the estimation of $H$ is 
 given by
	\begin{equation}
\widehat{H}=\frac{1}{2}\log _{2}\left( \frac{V_{n,a^{2}}}{V_{n,a}}\right),
\label{Hgorro}
\end{equation} 
where $a^{2}$ means the filter defined by $a^{2}=(a_0,0,a_1,0,a_2,0,...,0,a_k)$ of length $2k+1$ and order $L$ and 
$V_{n,a}:=\frac{1}{n}\sum_{i=0}^{n-k}\left( \sum_{j=0}^{k}a_{j}X_{i+j}\right)
^{2}$ is 
 the quadratic variation of the sample associated to a filter $a$.
 
Second, once we have estimated $H$ the  estimation of $\sigma$ is given by
 \begin{equation}
 \widehat{\sigma }=\left( \frac{-2V_{n,a}}{\Delta _{n}^{2\widehat{H}%
 	}\sum_{i=0}^{k}\sum_{j=0}^{k}a_{i}a_{j}\left\vert i-j\right\vert ^{2\widehat{%
 			H}}}\right) ^{1/2},  \label{sigmagorro}
 \end{equation}
 where $\Delta_n=\frac{T}{n}$.
 \\

\noindent \textbf{Step 2. Estimation of $\lambda_1$.}
 
If we define the function  
\begin{eqnarray*}
	U_{T}^{\left( n\right) }\left( \lambda_1,\lambda_2\right) &=&\frac{T}{n}%
	\sum_{i=1}^{n}h_{T}^{\left( n\right) }\left( iT/n,\lambda_1,\lambda_2  \right)
\end{eqnarray*}%
  in which $h_{T}^{\left( n\right) }$ is  

 \begin{equation*}
 h_{T}^{\left( n\right) }(x,\lambda_1,\lambda_2 )=\frac{1}{2\pi }\left( \log
 f^{\left( X\right) }\left( x,\lambda_1, \lambda_2 ,\widehat{\sigma},\widehat{H}\right) +\frac{I_{T}^{\left(
 		n\right) }(x)}{f^{\left( X\right) }\left( x,\lambda_1, \lambda_2 ,\widehat{\sigma},\widehat{H}\right) }%
 \right) w(x)
 \end{equation*}%
 where $f^{(X)}(x,\lambda_1,\lambda_2,\widehat{\sigma},\widehat{H})$ is the spectral density given in Section \ref{spectral_density} (evaluated in $H=\widehat{H}$ and $\sigma=\widehat{\sigma}$) and 
 \begin{equation*}
 I_{T}^{\left( n\right) }(x)=\frac{T}{2\pi }\left\vert \frac{1%
 }{n}\sum_{j=1}^{n}e^{\frac{ijTx}{n}}X_{\frac{jT}{n}}\right\vert ^{2}
 \end{equation*}%
 is the discretization of the periodogram. The weight function is $w(x)=\frac{|x|^c}{1+|x|^b}$ where $a$ and $b$ are parameters such that $c\geq 4$ and $b \geq c+3$.
 
Finally, the estimation of the parameters $(\lambda_1,\lambda_2 )$ are 
 \begin{eqnarray}
  (\widehat{\lambda }_1,\widehat{\lambda}_2)  &=&\arg \min_{(\lambda_1,\lambda_2 ) \in \Lambda
 }U_{T}^{\left( n\right) }\left( \lambda_1, \lambda_2\right) 
 \label{lambdasgorro}
 \end{eqnarray}%
 where $\Lambda \subset \left \{(\lambda_1, \lambda_2 ) \in \mathbb{R}^2 \ : \ 0<\lambda_1 < \lambda_2 \right \} $ some compact set.
 
 In   \cite{Kalemkerian_discreto} and \cite{chichi} can be found the proof of the asymptotic results concerning these estimators.

\subsection{Selection of a suitable value of $T$}
The value of $T$ gives us an idea about the unit measurement in which the observations  are taking. Although in 
every case it is natural to take a certain value of $T$ (e.g.,
if the observations are monthly and we have $120$ observations, it is natural to take $T=120$ months or $T=10$ y
ears) we can easily take any value of $T$ and interpret it in terms 
of the original time measure of the data. Therefore, we can take advantage of this fact and choose a
value of $T$ according to certain criteria. To be in accordance with the
asymptotic results 
we should choose a large value of $T$ such that $T/n$ is small.

Given a real time series, the decision of our test strongly depends on the choice of the parameter $T$. 
For example in Table \ref{arma11}, we have simulated different ARMA$(p,q)$ models 
where $p,q \leq 1$ and  we observe that (depending of the selection of $T$) we can
incorrectly reject the null hypothesis with an estimated probability higher or much higher than the 
significance level
($10 \%$ in Table \ref{arma11}).  The same problem is also present in
  other existing tests in the literature
(e.g., see   \cite{giraitis} and \cite{Teverovsky}). This happens because 
they also have a parameter (instead $T$) that must be chosen before carrying out the test. We 
show this fact in the Section \ref{section5}. 

\begin{figure}
	\centering
	\includegraphics[scale=0.7]{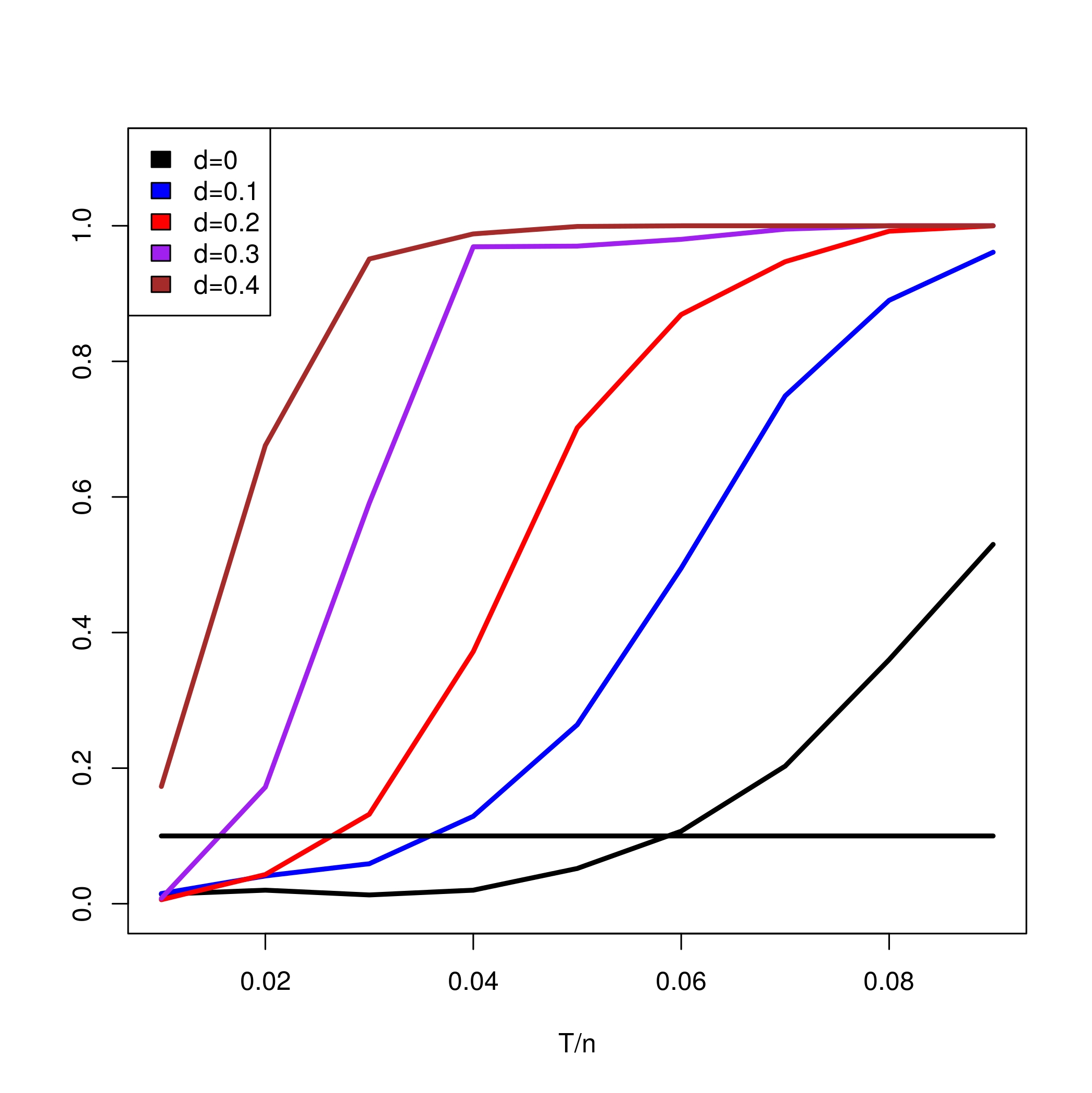} 
	\caption{Empirical power as a function of $T/n$  from $m=1000$ replications under $n=500$ observations of an ARFIMA$(1,d,1)$ model at significance level of $\alpha=0.1$ for different
	values of $d$, where the parameters are $\phi=0.8$ and $\theta=0.8.$ The horizontal line
	is the significance level of $\alpha=0.1$. 
	} 
	\label{size_fou}
\end{figure}

To deal with  this problem, we propose to select a value of $T$ such that the empirical probabilities of rejecting
the null hypothesis are less than or equal to the significance level under a 
certain family of short memory process.
Of course, it is impossible to consider every short memory process, hence we propose a value of $T$
such that the empirical probabilities of rejecting the null hypothesis under every ARMA$(p,q)$ process 
for $p,q \leq 1$ and $|\phi|, |\theta| \leq 0.8$ will be at most the significance level. 
The reason to consider $|\phi|, |\theta| \leq 0.8$ is based on the fact that 
the process is close to non-invertibility and non-causality in the other cases. 
In all of the simulated cases (both short memory and long memory processes), it is observed that the
empirical probabilities of rejecting the null hypothesis increases as $T$ increases, 
see  Figure 
\ref{size_fou}. Therefore,
to select the value  of $T$ to perform the test and to maximise the power of the test, we propose to select  a 
maximum value of $T$ such that the size 
of the test under every ARMA$(p,q)$ where $p,q \leq 1$ will be less than or equal to the significance level. The 
reason why 
 the ARMA$(p,q)$ models are included where $p,q \leq 1$  is due to the broad practical utility of these 
models but nevertheless it is a bit more general than other studies that only include 
AR$(1)$ models (see for instance \cite{Andrews}).

 Table \ref{arma11} shows that for sample size of $n=500$ we must consider $T \leq 0.06n$ as a possible 
 value of $T$ to carry out the test. It also shows us that the highest probabilities are obtained for the highest 
 value of
 the parameter $\phi$ and $\theta$. In Figure \ref{size_fou}, we show that for sample size of $n=500$ and
 $m=1000$ replications,  the empirical probabilities of rejecting 
 the null hypothesis under ARFIMA$(1,d,1)$ processes as 
 a function of $T$ when the AR and MA parameters are $0.8$ and $0.8$ respectively,
 for different values of $d$. The horizontal black line shows the signification level 
 ($\alpha=0.1$). It is observed that the power grows as the parameter $d$ grows (i.e., the test works well).
If $d=0$
 (short memory process), then the empirical probabilities of rejecting
 $H_0$ are greater than $0.1$ for values of $T > 0.06n$. Therefore, $T$ grows as the 
 empirical power grows and the optimal 
 value is reached at $T=0.06n$ (for $n=500$).
 Following this idea, in Table \ref{maximumT} we show the proposed value of $T$ to perform the test
for different signification levels and sample sizes. It is observed that the fraction $T/n$ is similar for 
every value of $n$ between $500$ and $5000$. 

\begin{table} 
	\centering
	\caption{Empirical probabilities, from $1000$ simulations, to reject the short range dependence at 
	significance level of $10\%$ under several ARMA$(p,q)$ alternatives where $p,q \leq 1$,
	 for $n=500$ (sample size) and different values of $T$. }
	
		\begin{tabular}{|c|cccc|}
		\hline
		$T$ & $0.1n$ & $0.05n$ & $0.01n$ &  $0.005n$  \\
		\hline
	AR$(1): \phi=0.2$	 & $0.000$   & $0.000$  & $0.000$  &$0.000$     \\
	AR$(1): \phi=0.5$	& $0.000$   & $0.000$  & $0.000$  & $0.000$    \\
	AR$(1): \phi=0.8$ & $0.001$  & $0.002$  & $0.000$  & $0.000$     \\
	MA$(1): \theta=0.2$	 & $0.000$   & $0.000$  & $0.000$  &  $0.000$    \\
	MA$(1): \theta=0.5$	 & $0.000$   & $0.000$  & $0.000$ & $0.000$     \\
	MA$(1): \theta=0.8$	 & $0.010$  & $0.062$  & $0.086$  & $0.000$     \\	
	ARMA$(1,1): (\phi, \theta)=(0.4,0.6)$	 & $0.000$   & $0.006$  & $0.067$  & $0.000$     \\
	ARMA$(1,1): (\phi, \theta)=(0.6,0.6)$	 & $0.288$  & $0.017$  & $0.024$  & $0.000$     \\
	ARMA$(1,1): (\phi, \theta)=(0.6,0.4)$	 & $0.006$   & $0.004$  & $0.039$ & $0.001$    \\
	ARMA$(1,1): (\phi, \theta)=(0.3,0.8)$	 & $0.000$  & $0.010$  & $0.059$ & $0.000$    \\
	ARMA$(1,1): (\phi, \theta)=(0.5,0.8)$	 & $0.163$  & $0.005$  & $0.016$ & $0.000$    \\
	ARMA$(1,1): (\phi, \theta)=(0.7,0.8)$	 & $0.987$   & $0.690$ & $0.033$  & $0.003$    \\
	ARMA$(1,1): (\phi, \theta)=(0.8,0.3)$	 & $0.397$ & $0.027$ & $0.022$  & $0.000$     \\
	ARMA$(1,1): (\phi, \theta)=(0.8,0.5)$	 & $0.943$  & $0.517$  & $0.002$  & $0.000$    \\
	ARMA$(1,1): (\phi, \theta)=(0.8,0.7)$	 & $1.000$   & $0.946$  & $0.127$  & $0.017$    \\	
		\hline

	\end{tabular}

\label{arma11}
\end{table}

\begin{table}
	\caption{Suitable value of the fraction $T/n$ for different sample sizes and significance levels.}
	\centering
	\begin{tabular}{|c|cccc|}
		\hline
		$ \alpha | n $ & $500$ & $1000$ & $3000$ & $5000$ \\
		\hline
		$0.010$ & $0.009$ & $0.045$ & $0.058$ & $0.062$ \\
		$0.025$ & $0.038$ & $0.050$ & $0.062$  & $0.065$ \\
		$0.050$ & $0.054$ & $0.056$  & $0.065$  & $0.068$  \\
		$0.075$ & $0.057$ & $0.060$  & $0.068$  & $0.069$  \\
		$0.100$ & $0.060$  & $0.063$ & $0.070$  & $0.072$ \\
		\hline
	\end{tabular}
	\label{maximumT}
\end{table}

\section{Comparison with other tests}\label{section5}
To compare the empirical power and the size of the proposed hypothesis test in this work 
(which we  call the \textit{FOU test})
with respect to other existing ones in the literature, we simulate different long and short memory processes
and we analyze the test performance. We start  by briefly describing each one of the tests that we  use to make 
the comparison.

\subsection{Other tests included in the comparison}
\textbf{1. The Lo test.}
This hypothesis test about short versus long range can be found
in \cite{Lo}. In this work, the author proposed the modified $R/S$ statistic, which is  defined  as
\begin{equation*}
Q_{n}(q):=\frac{1}{\widehat{\sigma }_{n}\left( q\right) }\left( \max_{1\leq
k\leq n}\sum_{i=1}^{k}\left( X_{i}-\overline{X}_{n}\right) -\min_{1\leq
k\leq n}\sum_{i=1}^{k}\left( X_{i}-\overline{X}_{n}\right) \right) , 
\end{equation*}%
where
\begin{equation*}
\widehat{\sigma }_{n}^{2}\left( q\right) =S_{n}^{2}+\frac{2}{n}%
\sum_{j=1}^{q}w_{j}\left( q\right) \sum_{i=j+1}^{n}\left( X_{i}-\overline{X}%
_{n}\right) \left( X_{i-j}-\overline{X}_{n}\right) \text{ being }w_{j}\left(
q\right) =1-\frac{j}{q+1}\text{, }q<n,
\end{equation*}%
and $S_n^{2}$ the sample variance.

It is important to note that $\widehat{\sigma }_{n}^{2}\left( q\right)
=S_{n}^{2}+2\sum_{j=1}^{q}w_{j}\left( q\right) \widehat{\gamma }_{j}$ where $%
\widehat{\gamma }_{j}$ the usual autocovariance estimator of lag $j$
and $q$ is the maximum lag to be considered. When $q=0$ we have the classical $R/S$ statistic first
used by Hurst (\cite{Hurst}).\\
\textbf{2. The V/S test.}
In  \cite{giraitis}, the author proposed a rescaled variance test based
on the $V/S$ statistic to test long range against short range. This  work shows
 that the statistic has a simpler asymptotic  than the modified rescaled range test
(Lo's test). The authors proposed the $%
M_{n}\left( q\right) $ statistic, which they called $V/S$ or rescaled
variance statistics and is defined by 
\begin{equation*}
M_{n}(q):=\frac{1}{\widehat{\sigma }_{n}^{2}(q)n^{2}}\left[
\sum_{k=1}^{n}\left( \sum_{j=1}^{k}\left( X_{j}-\overline{X}_{n}\right)
\right) ^{2}-\frac{1}{n}\left( \sum_{k=1}^{n}\sum_{j=1}^{k}\left( X_{j}-%
\overline{X}_{n}\right) \right) ^{2}\right] .
\end{equation*}%
The name $V/S$ comes from \textit{variance/S} because the statistic $M_{n}\left(
q\right) $ can be expresed as $M_{n}(q)=\frac{1}{\widehat{\sigma }%
_{n}^{2}(q)n}\widehat{\mathbb{V}}\left( S_{1}^{\ast },...,S_{n}^{\ast
}\right) $ where $S_{k}^{\ast }=\sum_{j=1}^{k}\left( X_{j}-\overline{X}%
_{n}\right) $ and $\widehat{\mathbb{V}}\left( S_{1}^{\ast },...,S_{n}^{\ast
}\right) =\frac{1}{n}\sum_{k=1}^{n}\left( S_{k}^{\ast }-\overline{%
S_{n}^{\ast }}\right) ^{2}.$ 

Observe that $Q_{n}(q)=\frac{1}{\widehat{\sigma 
}_{n}\left( q\right) }\left( \max_{1\leq k\leq n}S_{k}^{\ast }-\min_{1\leq
k\leq n}S_{k}^{\ast }\right) $, therefore the $M_{n}\left( q\right) $
statistic considers the sample variance of the values $S_{1}^{\ast
},...,S_{n}^{\ast }$ instead of the range of $S_{1}^{\ast },...,S_{n}^{\ast }$
adequately readjusted.  It is important to
highlight, that   this work  includes the asymptotic
distribution under short memory and long memory dependence for the
statistics $M_{n}\left( q\right) $ and $Q_{n}\left(
q\right) .$
Under some conditions of $q$ and the fourth cumulants 
the authors proved that $M_{n}$ converges to $F_{KS}\left( \sqrt{\pi 
}x\right) $, where $F_{KS}$ is the asymptotic distribution of the
Kolmogorov-Smirnov statistic.
\\
\textbf{3. The Gromykov et al test.}
In \cite{Gromykov}, the authors proposed to split the original sample $%
X_{1},X_{2},...,X_{n}$ into $m$ blocks of size $l$ and construct the
periodogram of the entire sample given by  $I_{n}\left( \lambda \right) =%
\frac{1}{2\pi n}\left\vert \sum_{j=1}^{n}X_{j}e^{ij\lambda }\right\vert ^{2}$
and the periodogram of each block $I_{n,i}$ for $i=1,2,...,m$ to then work
with 
\begin{equation*}
Q_{n,m}(s)=\sum_{j=1}^{s}\frac{I_{n}\left( \lambda _{j}\right) }{\frac{1}{m}%
\sum_{i=1}^{m}I_{n,i}\left( \lambda _{j}\right) }
\end{equation*}%
as a test statistic, where $s$ is the number of Fourier frequencies to
consider and $\lambda _{j}=\frac{2\pi j}{n}$.

If $\left\{ X_{t}\right\} $ is a stationary linear process that is defined as  $%
X_{t}=\sum_{j=0}^{+\infty }a_{j}\varepsilon _{t-j}$ where $\left\{
\varepsilon _{j}\right\} $ are i.i.d. with zero mean and variance $\sigma
^{2}>0,$ is a short memory process in the sense that it has an spectral
density of the form $f\left( \lambda \right) =\left\vert \lambda \right\vert
^{-2d}g\left( \lambda \right) $ \ for $d =0$ and $\sum_{j=1}^{+\infty
}\left\vert a_{j}\right\vert <+\infty $, then  for fixed $s$ the authors proved that  $%
Q_{n,m}(s)\overset{d}{\rightarrow }Q\left( s\right) $  when  $m\rightarrow
+\infty $ being $m=o(n)$ where $Q\left( s\right) $ has Gamma distribution
with parameters $\left( s,1\right) .$ The details and the
asymptotic results for the case $0<d<1/2$ can be found in \cite{giraitis}.
According with the suggestion in \cite{Gromykov}, we use $m=\sqrt{n}$ as a number of blocks. From now on,
we will call the
Gromykov et al test as the $Q$ test. \\
 \textbf{4. The Lobato--Robinson test.}
 The Lobato-Robinson test is based on the statistic defined by LR$=t^{2}$
where 
\begin{equation*}
t=\sqrt{m}\frac{\sum_{j=1}^{m}\nu _{j}I\left( \lambda _{j}\right) }{%
\sum_{j=1}^{m}I\left( \lambda _{j}\right) }
\end{equation*}%
where $\upsilon _{j}=\log
j-\frac{1}{m}\sum_{i=1}^{m}\log i $, $\lambda _{j}=\frac{2\pi j}{n}$ and $I\left( \lambda \right) $ is the
periodogram. 
Under certain conditions of the number of blocks $ m$, the
authors proved that if the spectral density is twice bounded differentiable
near $0$ and $d=0$, then the value $t$ converges in distribution of an standard normal
distribution.

\begin{figure}[H]
	\centering
	\includegraphics[scale=0.5]{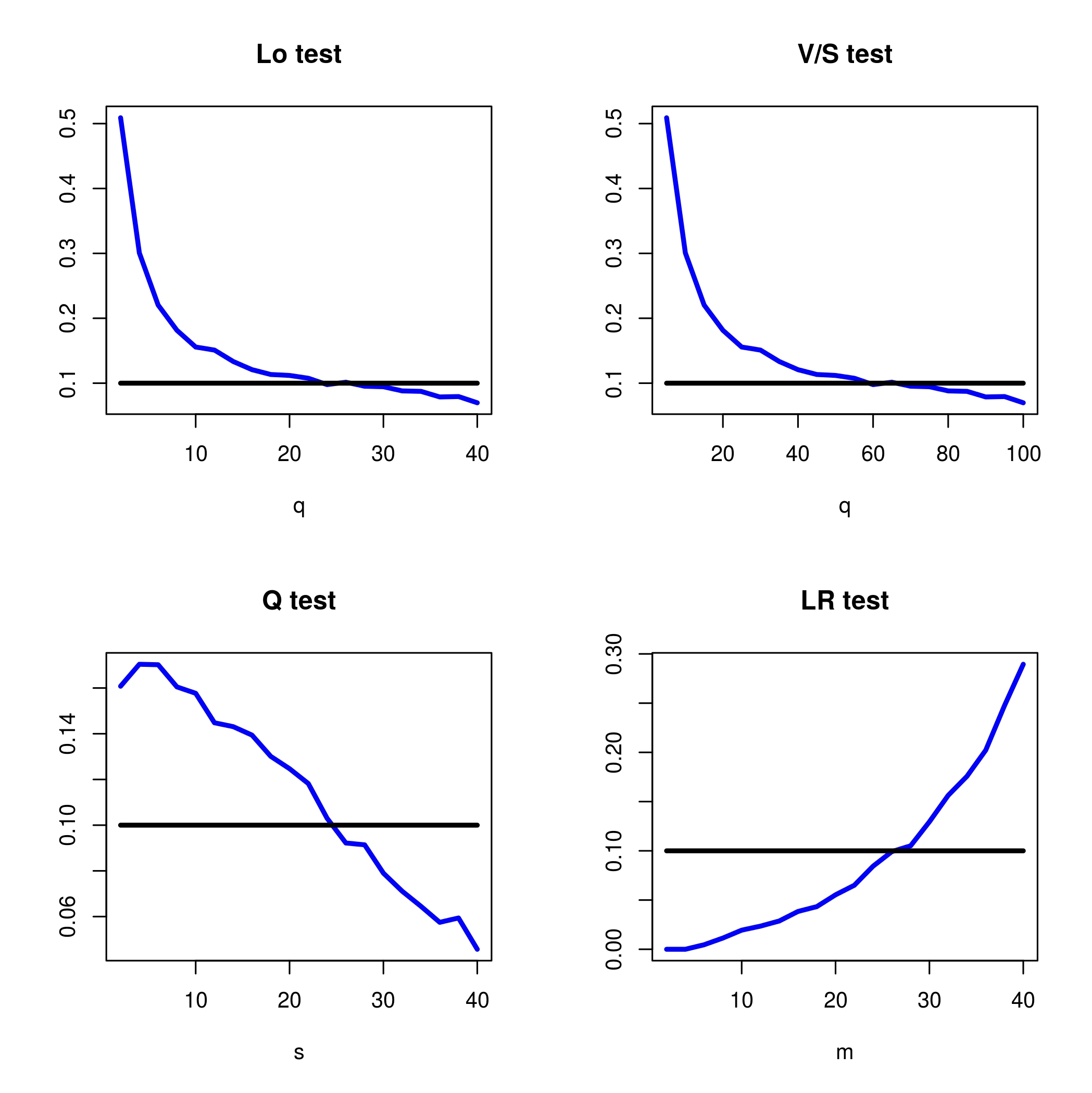} 
	\caption{Empirical probabilities for $m=1000$ replications of reject the null hypothesis for an ARMA$(1,1)$ model with parameters $\phi=\theta=0.8$ at 
	the $10\%$ level of significance where the sample size is $n=1000$.
	} 
	\label{size_other_tests}
\end{figure}

\begin{remark}
 It is important to note that all of the considered tests have a parameter to select before applying it ($q$ in Lo and V/S, $m$ and $s$ in $Q$ test
 and $m$ in LR). The correct choice of the parameter is very difficult because if we select a value that is 
 too small
 or too high, then  each of the considered tests
 will have a bias to incorrectly reject the null hypothesis in some cases or incorrectly non-reject the null hypothesis in others. 
\label{decreasing}
\end{remark}
To illustrate the affirmation given  in Remark \ref{decreasing}, we observe in Figure \ref{size_other_tests} the
behavior of the empirical probability of rejecting the null hypothesis  in the case in which the observed time 
series corresponds to an ARMA$(1,1)$ 
model where $\theta=\phi=0.8$. The empirical probability of rejecting the null hypothesis decreases as $q$
increases in Lo and V/S tests, decreases as $s$ increases in the Gromykov et al test and 
increases as $m$ increases in the LR test. The same  behavior is repeated for any other simulated
time series (short or long memory).

\subsection{Short and long range processes considered in the comparison}
To carry out a  comparison between the FOU test and the other tests described in this subsection, we  use the
optimum value of $T$ according to Table \ref{maximumT}. 
 For the other tests (Lo, V/S, $Q$ and LR) we  use the same criterion proposed in Subsection 4.2 and given in Table
 \ref{q_optimal}.
 In this way, all of the  tests that we have considered have an empirical probability  less than or equal to
 $0.1$ of incorrectly
 rejecting the null hypothesis under the same family of short 
 range dependence models, that is, under every ARMA$(1,1)$ where $|\phi|, |\theta| \leq 0.8$.
We have made the comparison looking the power under several long range alternatives and the size of
 the considered tests under several short range processes. 
 In the comparison, we have considered  ARFIMA, LARCH, FGN and FOU models.
 All these processes include (depending on 
 the value of their parameters) both short and long range memory,
 whose definitions are outlined below.\begin{enumerate}

                                     \item A LARCH model is defined by $X_{k}=r_{k}^{2},$ where $%
r_{k}=\sigma _{k}\varepsilon _{k}$, $\sigma _{k}^{2}=\left( \alpha
+\sum_{j=1}^{+\infty }\beta _{j}r_{k-j}\right) ^{2}$ and $\left\{
\varepsilon _{k}\right\} $ are white noise.  Depending on the speed at which the coefficients tend to zero
we can have a short or long memory process.
\item  $\left\{ X_{t}\right\} $ is an ARFIMA$\left( p,d,q\right) $
process when $\left( 1-B\right) ^{-d}X_{t}$ is an ARMA$\left( p,q\right) $
process. If $d=0$  we have an ARMA$\left( p,q\right) $ short memory process
and when $0<d<1/2$ we have a long memory process.
\item The fractional Gaussian process (FGN$\left( H\right) $) is defined as a
stationary Gaussian centered process $\left\{ X_{t}\right\} _{t\in \mathbb{N}%
}$ such that the autocovariance function is given by $\mathbb{E}\left(
X_{t}X_{0}\right) =\frac{1}{2}\left( \left( t+1\right) ^{2H}-2t^{2H}+\left(
t-1\right) ^{2H}\right) .$ It is known that when $H>1/2$,  we have a
long memory process and if $H\leq 1/2$ then we have a short memory process.
\item The fractional Ornstein--Uhlenbeck process FOU$(\lambda,\sigma,H)$ was defined in Section \ref{section2}.
When $H \leq 1/2$ the process has short memory and when $H>1/2$ long memory. When $H=1/2$ we have 
the classical Ornstein--Uhlenbeck process.

                                                                                                         \end{enumerate}

 \subsection{Power comparison}
 
   \begin{table}
  \caption{Optimal value of $q$ for the Lo test and V/S test at the significance level of
  $10 \%$ ($ 5 \%$), optimal value of $s$ for the $Q$ test, and optimal value of $m$ for the LR test, 
  in function of sample size $n$.}
  \centering
  \begin{tabular}{|c|cccc|}
  \hline $n$ & $500$ & $1000$ & $3000$ & $5000$ \\
  \hline Lo & $16 (15)$ & $ 23 (22)$ & $38 (37)$ & $51 (50)$ \\
  V/S & $ 44 (32)$ & $66 (47)$ & $90 (80)$ & $ 128 (100)$ \\
  $Q$ & $15 (7)$ & $25 (9)$ & $21 (11)$ &  $22 (12)$ \\
  LR & $17 (17)$ & $27 (28)$ & $57 (54)$ & $81 (77)$\\
  \hline
  \end{tabular}               
\label{q_optimal}
 \end{table}
 In Table \ref{power_arfima} we show the power at $10 \%$ of the tests 
 under different ARFIMA$(1,d,1)$ alternatives for different values of $d>0$ (i.e., long memory processes).
 The power of each test was obtained from $1000$ replications. Table \ref{power_arfima} shows that the FOU test
 detects long memory when none of the parameters is close to zero. 
 It is also observed that the performance of the test improves as both parameters increase their value and when 
 $d$ increases (as expected).
 When the value of both parameters is high, the FOU test has the best performance, while in
 almost all of the other cases, the LR test obtain the best results.

\begin{table}
	\caption{Power comparison from $1000$ simulations under several ARFIMA$(1,d,1)$
		 alternatives and different sample sizes ($n$) at level $10 \%$. All calculations 
		 were performed using the values of the test parameters suggested in \ref{maximumT} and \ref{q_optimal}. }

\begin{small}\begin{tabular}{|c|c||c|c|c||c|c|c||c|c|c|}
\hline
\multicolumn{11}{|c|}{\ \ \ \ \ \ \ \ \ \ \ \ \ \ \ \ \ \     $(0.3,0.8)$ \ \ \ \ \ \ \ \ \ \ \ \ \   \ \ \ \ \ \ \  \ $(0.5,0.8)$ \ \ \ \ \ \ \ \ \ \  \ \ \ \ \ \ \ \ \ \  \ \ $(0.7,0.8)$ } \\
\hline
$d$ & $n$ & $500$ & $1000$ & $5000$ &  $500$ & $1000$ & $5000$ & $500$ & $1000$ & $5000$ \\
\hline
$0.1$ & FOU & $0.078$  & $0.024$ & $0.005$  & $0.016$ & $0.021$ & $0.000$  & $0.007$   & $0.048$ & $0.020$  \\
& V/S & $\textbf{0.155}$  & $0.175$ & $0.270$ & $\textbf{0.189}$  & $\textbf{0.207}$ & $0.294$  & $\textbf{0.181}$  & $0.209$ & $0.277$ \\
& Lo & $0.144$ & $\textbf{0.194}$ & $0.317$ & $0.165$ & $0.201$ & $0.319$  & $0.171$  & $0.191$  & $0.317$ \\
& Q & $0.054$ & $0.055$ & $0.206$ & $0.059$ & $0.066$ & $0.263$ & $0.088$ & $0.075$ & $0.273$ \\
& LR & $0.109$ & $0.192$ & $\textbf{0.489}$  & $0.114$ & $0.205$ & $\textbf{0.540}$ & $0.142$ & $\textbf{0.254}$ & $\textbf{0.542}$ \\

\hline
$0.2$ & FOU & $0.079$   & $0.029$ & $0.020$  & $0.015$ & $0.020$ & $0.000$  & $0.291$  & $0.318$ & $0.330$  \\
& V/S & $0.259$  & $0.333$ & $0.496$  & $0.284$  & $0.337$ & $0.472$ & $0.310$  & $0.209$ & $0.452$ \\
& Lo & $0.280$ & $0.358$ & $0.608$  & $0.257$ & $0.358$ & $0.585$ & $0.303$  & $0.392$  & $0.587$ \\
 & Q & $0.062$ & $0.068$ & $0.587$ & $0.050$ & $0.073$ & $0.586$ & $0.108$ & $0.092$ & $0.615$ \\ 
& LR & $\textbf{0.300}$ & $\textbf{0.490}$ & $\textbf{0.940}$ & $\textbf{0.301}$ & $\textbf{0.521}$ & $\textbf{0.929}$ & $\textbf{0.330}$ & $\textbf{0.559}$ & $\textbf{0.925}$ \\

\hline 
$0.3$ & FOU & $0.082$  & $0.038$ & $0.015$  & $0.071$ & $0.037$ &  $0.018$ & $\textbf{0.723}$  & $\textbf{0.780}$ & $0.895$\\
& V/S & $0.390$  & $0.464$ & $0.689$  & $0.396$ & $0.493$ & $0.667$ & $0.394$  & $0.487$ & $0.693$ \\
& Lo & $0.408$ & $0.555$ & $0.796$  & $0.448$ & $0.577$ & $0.818$  & $0.438$ & $0.564$ & $0.808$  \\
& Q & $0.047$ & $0.060$ & $0.886$ & $0.070$ & $0.071$ & $0.883$ & $0.077$ & $0.077$ & $0.889$ \\
& LR & $\textbf{0.504}$ & $\textbf{0.748}$ & $\textbf{0.999}$ & $\textbf{0.540}$ & $\textbf{0.767}$ & $\textbf{0.998}$ & $0.592$ & $0.774$ & $\textbf{0.998}$ \\

\hline
$0.4$ & FOU & $0.057$  & $0.031$ & $0.005$  & $0.324$ & $0.350$ & $0.659$  & $\textbf{0.987}$   & $\textbf{0.998}$ & \textbf{1.000}   \\
& V/S & $0.502$  & $0.603$ & $0.798$ & $0.510$ & $0.611$ & $0.798$ & $0.509$  & $0.635$ & $0.783$ \\
& Lo & $0.525$  & $0.688$ & $0.907$ & $0.551$ & $0.698$ & $0.924$  & $0.571$ & $0.678$ & $0.903$  \\
& Q & $0.060$ & $0.044$ & $0.982$ & $0.051$ & $0.044$ & $0.989$ & $0.043$ & $0.047$ & $0.985$ \\
& LR & $\textbf{0.702}$ & $\textbf{0.923}$ & $\textbf{1.000}$ & $\textbf{0.732}$  & $\textbf{0.927}$ & $\textbf{1.000}$ & $0.770$ & $0.946$ & $0.999$ \\
\hline

\end{tabular}             \end{small}
\\
\begin{small}
\centering
\begin{tabular}{|c|c||c|c|c||c|c|c||c|c|c|}
	\hline
	\multicolumn{11}{|c|}{\ \ \ \ \ \ \ \ \ \ \ \ \ \ \ \ \ \     $(0.8,0.3)$ \ \ \ \ \ \ \ \ \ \ \ \ \   \ \ \ \ \ \ \  \ $(0.8,0.5)$ \ \ \ \ \ \ \ \ \ \  \ \ \ \ \ \ \ \ \ \  \ \ $(0.8,0.7)$ } \\
	\hline
	$d$ & $n$ & $500$ & $1000$ & $5000$ &  $500$ & $1000$ & $5000$ & $500$ & $1000$ & $5000$ \\
	\hline
	$0.1$ & FOU & $0.016$   & $0.008$  & $0.010$   & $0.081$ & $0.036$ & $0.061$  & $\textbf{0.426}$  & $\textbf{0.504}$ & $\textbf{0.934}$ \\
	& V/S & $0.200$ & $0.231$ & $0.286$ & $0.228$ & $0.204$ & $0.290$  & $0.198$ & $0.239$ & $0.282$ \\
	& Lo & $0.181$ & $0.219$ & $0.331$ & $0.203$ & $0.267$ & $0.343$ & $0.196$ & $0.206$ & $0.343$ \\
	& Q & $0.116$ & $0.117$ & $0.304$ & $0.119$ & $0.143$ & $0.327$ & $0.120$ & $0.146$ & $0.324$ \\
	& LR & $\textbf{0.269}$ & $\textbf{0.339}$ & $\textbf{0.619}$ & $\textbf{0.250}$ & $\textbf{0.342}$ & $\textbf{0.626}$ & $0.257$ & $0.356$ & $0.617$ \\
	
	\hline
	$0.2$ & FOU & $0.182$  & $0.126$  & $0.312$  & $\textbf{0.705}$ & $\textbf{0.808}$ & $0.841$  & $\textbf{0.876}$  & $\textbf{0.933}$ & $\textbf{0.993}$ \\
	& V/S & $0.340$  & $0.343$ & $0.484$ & $0.341$  & $0.351$ & $0.476$  & $0.305$  & $0.351$ & $0.491$ \\
	& Lo & $0.335$ & $0.421$ & $0.600$  & $0.337$ &$0.425$ & $0.625$  & $0.321$  & $0.423$ & $0.593$ \\
	& Q & $0.088$ & $0.123$ & $0.686$ & $0.107$ & $0.138$ & $0.654$ & $0.100$ & $0.121$ & $0.670$ \\
		& LR &$\textbf{0.462}$ & $\textbf{0.633}$ & $\textbf{0.967}$ & $0.460$ & $0.645$ & $\textbf{0.962}$ & $0.452$ & $0.644$ & $0.970$ \\
	\hline 
	$0.3$ & FOU & $\textbf{0.795}$  & $\textbf{0.870}$ & $0.983$  & $\textbf{0.976}$ & $\textbf{0.997}$ & $\textbf{0.999}$  & $\textbf{0.991}$  & $\textbf{1.000}$ & $\textbf{1.000}$ \\
& V/S & $0.416$ & $0.498$ & $0.676$ & $0.575$ & $0.482$ & $0.644$  & $0.421$ & $0.499$ & $0.697$ \\
	& Lo & $0.474$  & $0.575$ & $0.808$ &   $0.428$ & $0.594$ & $0.828$ & $0.484$ & $0.579$ & $0.801$ \\
	& Q & $0.073$ & $0.108$ & $0.914$ & $0.061$ & $0.094$ & $0.925$ & $0.075$ & $0.110$ & $0.915$ \\
		& LR & $0.679$ & $0.844$ & $\textbf{1.000}$ & $0.659$ & $0.832$ & $0.998$ & $0.667$ & $0.843$ & $0.997$ \\
	\hline
	$0.4$ & FOU & $\textbf{0.995}$  & $\textbf{1.000}$ &   $\textbf{1.000}$ & $\textbf{1.000}$ & $\textbf{1.000}$ & $\textbf{1.000}$  & $\textbf{1.000}$  & $\textbf{1.000}$  & $\textbf{1.000}$ \\
	& V/S & $0.531$  & $0.612$& $0.800$  & $0.420$  & $0.579$ & $0.796$  & $0.526$  & $0.602$ & $0.806$\\
	& Lo & $0.597$ & $0.709$ & $0.928$ & $0.618$  & $0.722$ & $0.908$ & $0.595$ & $0.699$ & $0.922$  \\
	& Q & $0.055$ & $0.048$ & $0.987$ & $0.049$ & $0.055$ & $0.990$ & $0.044$ & $0.062$ & $0.993$\\
		& LR & $0.789$ & $0.940$ & $\textbf{1.000}$ & $0.786$ & $0.952$ & $\textbf{1.000}$ & $0.856$ & $0.953$ & $\textbf{1.000}$ \\
	\hline

\end{tabular}             \end{small}
\label{power_arfima}
\end{table}

In Table \ref{power_fgm_larch} we have considered  FGN$(H)$ for different values of $H>1/2$ and the LARCH$\left( 0,d,0\right) $ long memory process
where $\beta _{j}=j^{d-1}$ and $\alpha =0.1.$ For these families, we do not include the performance of the FOU
test because the estimation of $H$ is clearly less 
than $1/2$, and therefore the test non-rejects the null hypothesis. This occurs due
to what was observed in Remark \ref{H_interpretation} and Remark \ref{H_reject}. The
LR test obtains the best performance in all cases for a sample size of $n=5000$.
\begin{table}
	\caption{Power comparison from $1000$ simulations under several fractional Gaussian noise (FGN) and LARCH$(0,d,0)$
		alternatives and different sample sizes ($n$) at level $10 \%$. All calculations 
		were performed using the values of the test parameters suggested in \ref{maximumT} and \ref{q_optimal}. }
	\centering
	\begin{tabular}{|c|c||c|c|c|||c|c|c|c|c|}
		\hline
		\multicolumn{10}{|c|}{\ \ \ \ \ \ \ \ \ \ \ \ FGN$(H)$ \ \ \ \ \ \ \ \ \ \ \ \ \ \ \ \ \ \ \ \ \ \ \ \ \  \ LARCH$(0,d,0)$ \ \ \  \ \ \ \ \ \  } \\
		\hline
		$H$ & $n$ & $500$ & $1000$ & $5000$ & $d$ & $n$ & $500$ & $1000$ & $5000$  \\
		\hline

		$0.6$ & V/S & $0.122$  & $0.154$  & $0.248$  & $0.1$ & V/S  & $0.044$ & $0.066$ & $0.099$     \\
	
	     & Lo & $\textbf{0.147}$  & $\textbf{0.208}$  & $0.312$  &  & Lo & $0.013$  & $0.029$  & $0.040$    \\
	     
	      & Q & $0.020$  & $0.049$  & $0.200$  &  & Q & $\textbf{0.126}$  & $\textbf{0.146}$  & $0.212$     \\
	
		 & LR & $0.099$  & $0.180$  & $\textbf{0.487}$  &  & LR & $0.075$  & $0.122$  & $\textbf{0.271}$     \\
		\hline

	 $0.7$ & V/S &  $0.210$ & $0.266$ & $0.423$  & $0.2$ & V/S & $0.041$ & $0.063$ & $0.107$     \\
	
	& Lo & $0.273$   & $0.345$  & $0.602$  &  & Lo & $0.009$  & $0.033$  & $0.062$     \\
	
	& Q & $0.050$  & $0.054$  & $0.552$  &  & Q & $\textbf{0.138}$  & $\textbf{0.158}$  & $0.207$     \\
	
	& LR & $\textbf{0.281}$  & $\textbf{0.450}$  & $\textbf{0.913}$  &  & LR & $0.084$  & $0.153$ & $\textbf{0.354}$     \\
		
		\hline

	$08$ & V/S & $0.285$  & $0.381$  & $0.640$  & $0.3$ & V/S & $0.053$  & $0.058$  & $0.151$     \\
	
	& Lo & $0.408$  & $0.541$  & $0.776$  &  & Lo & $0.017$  & $0.034$  & $0.073$     \\
	
	& Q & $0.051$   & $0.057$  & $0.866$ &  & Q & $\textbf{0.141}$  & $0.151$  & $0.215$     \\
	
	& LR & $\textbf{0.494}$  & $\textbf{0.766}$  & $\textbf{0.999}$  &  &  LR & $0.088$  & $\textbf{0.152}$ & $\textbf{0.381}$    \\
		\hline

 $0.9$	& V/S &  $0.387$  & $0.494$ & $0.736$  & $0.4$ & V/S & $0.050$  & $0.058$  & $0.116$     \\
	
	& Lo & $0.532$  & $0.685$  & $0.903$  &  & Lo & $0.026$  & $0.031$ & $0.066$     \\
	
	& Q & $0.041$   & $0.042$  & $0.980$ &  & Q & $0.144$  & $0.148$  & $0.242$     \\
	
	& LR & $\textbf{0.721}$   & $\textbf{0.910}$  & $\textbf{1.000}$  &  & LR & $\textbf{0.095}$  & $\textbf{0.186}$  & $\textbf{0.383}$    \\
		\hline

	\end{tabular}
	
	\label{power_fgm_larch}
\end{table}

\subsection{Size comparison}
In this subsection, we have considered the  LARCH$%
\left( 1,0,1\right) $ short memory process, where $\phi =0.1,$ $\theta =0.2,$
$\alpha =0.1,$ $\beta _{j}=\phi ^{j-1}\left( \phi -\theta \right) $ and the $\varepsilon_{i}$ are i.i.d. normal standard variables, FGN$(H=0.5)$ and 
Ornstein--Uhlenbeck process where $\sigma=1$ and $\lambda=0.8$. Concerning  the ARMA process,
we know that  none of the hypothesis tests that we have considered fail to reject the null hypothesis for every
ARMA$(1,1)$ (in the sense that the percentage of reject the hypothesis of short memory process is not 
greater that the significance level) where $|\theta|, |\phi| \leq 0.8$. 
For this reason we have considered ARMA$(p,q)$ where $p$ or $q$ are greater than $1$.
We can draw important conclusions from Table \ref{size}. The V/S, Lo, Q and LR test work well under 
FGN, but they are terribly wrong in some cases of AR(2) alternatives and 
under Ornstein--Uhlenbeck processes, and the error get worse as the sample size increases. The FOU test
is never wrong under the FGN and LARCH short memory models. Under the ARMA models the FOU test
work well. Only under the Ornstein--Uhlenbeck processes observed in $[0,T]$  the empirical percentage 
of rejection of the null hypothesis in the FOU test is sligtlhy greater
than the size. In addition, this percentage increases as $T$ decreases. This is to be expected because if $T
\rightarrow 0$, then the autocovariances of the process goes to
$\mathbb{V}(X_0)>0$, and therefore we get closer to a long memory process.
\begin{table}
	\caption{Size at level $10 \%$ of any considered test under several short range models. The empirical
	probabilities were calculated from $1000$ replications. The ARMA$(0.4,0.55)$ means an AR$(2)$ where
	$\phi=(0.4,0.55)$.
	The OU case means an Ornstein--Uhlenbeck process where $\sigma=1$ and $\lambda=0.8$, observed 
	in $[0,T]$ being $T=100,50$ and $10$.
	All of the calculations 
		were performed using the values of the test parameters suggested in \ref{maximumT} and \ref{q_optimal}. }
	\centering
{\scriptsize 	\begin{tabular}{|c||c|c||c|c||c|c||c|c||c|c|}
		\hline
		\multicolumn{11}{|c|}{ \ \ \ \ \ \ \ \ \ \ \ \ \ \ \ \ \ \   \ \ \ \ \ FOU \ \ \ \ \ \ \ \  \ \ \ \ \ \ \ \  V/S \ \ \ \ \ \ \ \ \ \ \ \ \  \ \  \ \ \ \  Lo \ \ \ \ \ \ \ \ \ \ \ \ \ \ \ \ \ \ \ Q \ \ \ \  \ \ \ \ \ \ \ \ \ \ \ \ \ \ \ \ LR } \\
		\hline
		  & $1000$ & $5000$  & $1000$ & $5000$  & $1000$ & $5000$ & $1000$ & $5000$ & $1000$ & $5000$ \\
		\hline
		AR$(0.4,0.55)$ & $0.000$  & $0.000$  & $0.300$ & $0.214$   & $0.542$  & $0.357$  & $0.003$  & $0.892$    & $0.988$ & $1.000$   \\
		FGN$(H=0.5)$ & $0.000$ &  $0.000$ & $0.075$  & $0.120$   & $0.056$  & $0.084$  & $0.015$  & $0.032$    & $0.046$  & $0.079$   \\
		OU$(T=100)$ & $0.022$  & $0.125$   & $0.125$  & $0.311$  & $0.264$  & $0.680$ & $0.660$ & $1.000$   & $0.127$  & $0.974$   \\
		OU$(T=50)$ & $0.053$ &  $0.168$ & $0.230$ & $0.576$   & $0.481$ & $0.898$ & $0.026$ & $0.952$   & $0.961$  & $1.000$   \\
		OU$(T=10)$ & $0.123$ & $0.175$  &  $ 0.697$ & $0.959$  & $0.915$ & $0.997$ & $0.001$  & $0.642$    & $0.999$  & $1.000$   \\
		LARCH$(1,0,1)$ & $0.000$ & $0.000$  & $0.066$ & $0.099$   & $0.029$  & $0.040$ & $0.146$ & $0.212$    & $0.044$  & $0.060$   \\

		\hline

	\end{tabular}}
	\label{size}
\end{table}

\section{Application to  real data} \label{section6}
\begin{figure}[H]
	\centering
	\includegraphics[scale=0.4]{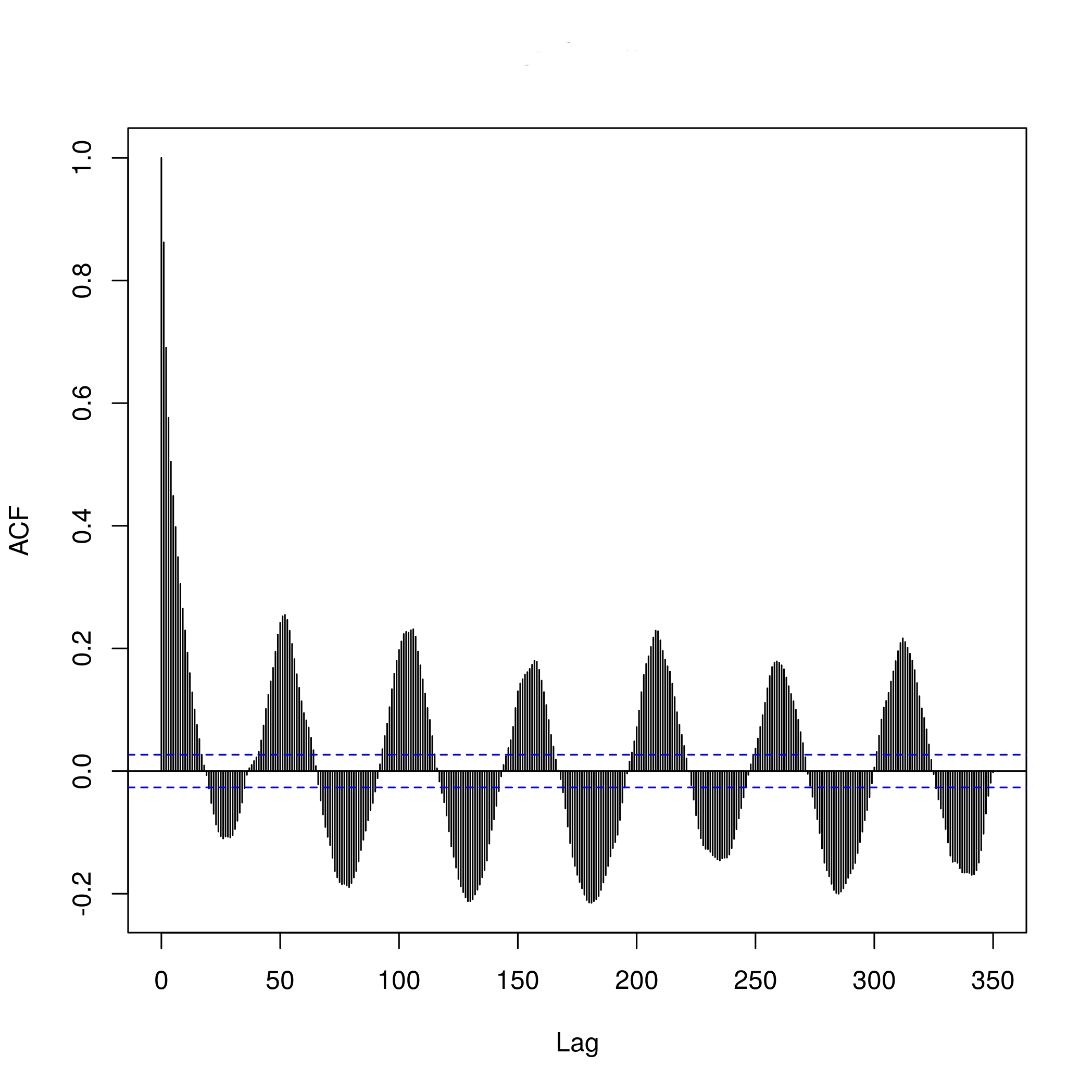} 
		\includegraphics[scale=0.4]{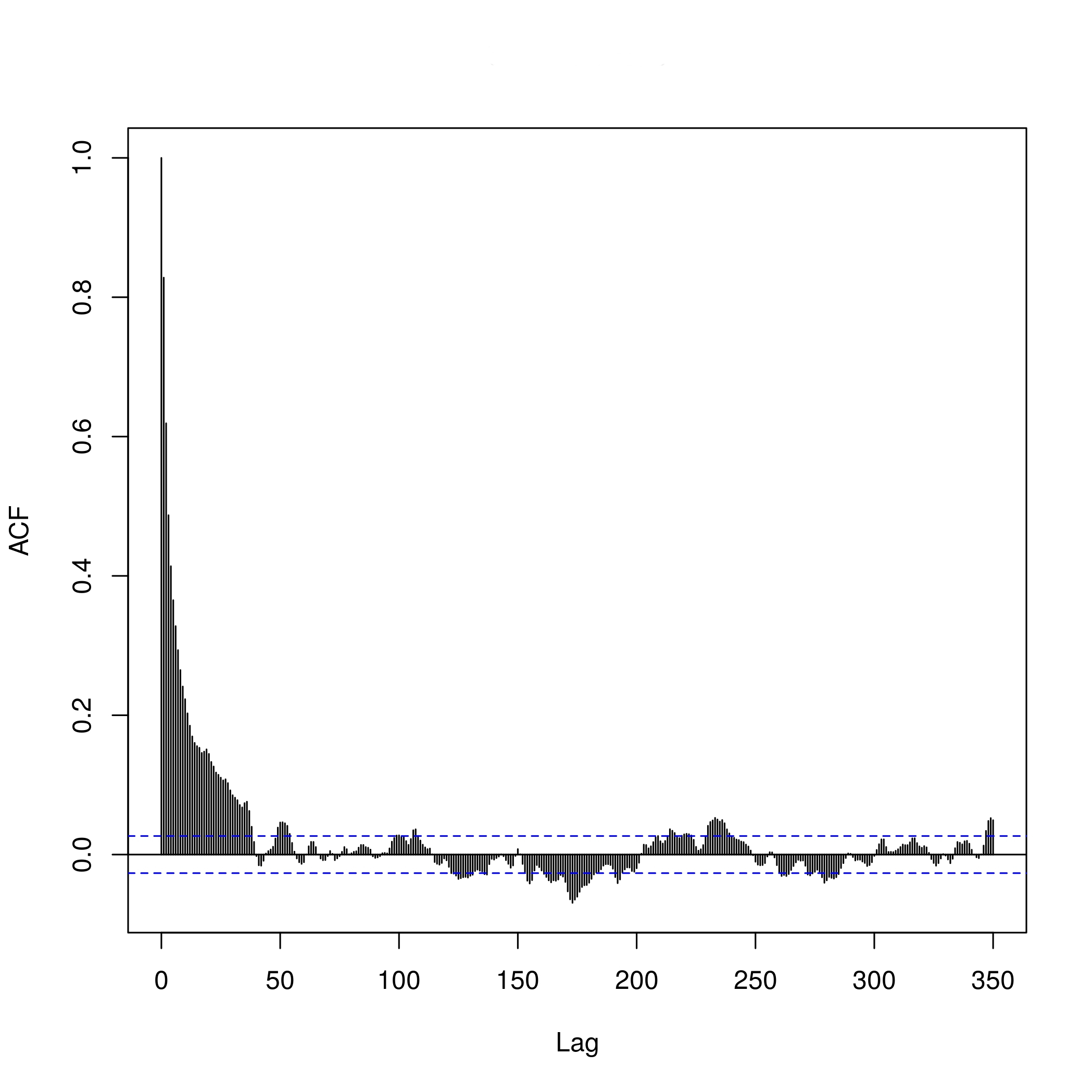} 

	\caption{Autocorrelation of the original time series (left) and seasonally adjusted time series (right).} 
	\label{acf_affluent}
\end{figure}

In this section  we analyze an empirical application  of a  time series that has already been studied and modeled
correctly with an ARFIMA long memory model in the work \cite{affluent}. The dataset consists of weekly
measurements of affluent energy generated by hydroelectric dams
in Uruguay between the first week of 1909 and the last week of 2012. The time series has length $5408$.
Each observation corresponds to the  weekly inflow  energy generated by the three Uruguayan
dams, measures in MWh. This time series is strongly related with the time series generated by the dam 
contributions and is also a good fit to a large memory 
model (\cite{represas}). After being seasonally adjusted  and centered, this time series 
 has a good fit to an ARFIMA$(3,d,1)$ long memory process (see \cite{affluent} for details).
 In Figure \ref{acf_affluent} we show the autocorrelation function of this time series before and after
 adjusted seasonally.
 
 The main improvement of modelling this time series using long memory processes can be seen in 
 the intensity curves (see the definition in \cite{represas}). The intensity curves can be used as a measure of persistence
 of droughts, especially the minimum curve, essential for energy planning. In
 \cite{Graneri} (Figure 14 and Figure 55) the intensity curves modeled by ARMA
 processes are shown, while in \cite{represas} (Figure 20) they are obtained by ARFIMA
 processes.
  
To perform the proposed test we have used (according to Table \ref{maximumT})
$T=0.072n=389.376.$ The FOU test clearly rejects
the hypothesis of short memory dependence. In Table \ref{affluent_test}, we show the parameter estimation for this 
real time series fitted to a FOU$(\lambda_1,\lambda_2,\sigma,H)$ model. The estimated model is a 
FOU$(\lambda_1=0,\lambda_2=5.8311,\sigma =4.4593,H=0.7078)$ model, 
which corresponds to a long memory process. We also observe that the estimated value of $H=0.7078$ is 
much higher than the critical value of the test. In this empirical application, 
the hypothesis test result is consistent with the previous results in \cite{affluent}.
\begin{table}
 \caption{Parameter estimations for the affluent energy data series fitted to   a FOU$(\lambda_1,\lambda_2,\sigma,H)$ model observed in $[0,T]$ for $T=389.376$.}
 \centering
 \begin{tabular}{|c|c|c|c||||c|c|}
 \hline $\hat{\lambda_1}$ & $\hat{\lambda_2}$ & $\hat{\sigma}$ & $\hat{H}$ & critical value at $10\% $ & p-value \\
 \hline $10^{-11}$ & $5.8311$ & $4.4593$ & $0.7078$ & $0.5033$ & $0.0000$\\
 \hline
  
 \end{tabular}
\label{affluent_test}
\end{table}

\section{Conclusions} \label{section7}

In this work we present a new hypothesis test to contrast short memory versus long memory in time series, which is based on  the  Fractional Iterated Ornstein--Uhlenbeck processes.
We  present an implementation of the test and we  carry out a simulation study that includes several families of processes of both short memory and long memory. We also compare the
results with other hypothesis tests. In addition, we propose a suitable value of $q$ to be used for the Lo test and V/S test, $m,s$ for the Q test and $m$ for LR test. 
With this election of the parameters, all of the considered tests maintain empirical probabilities of rejecting $H_0$ under every ARMA$(1,1)$ where $|\phi|, |\theta| \leq 0.8$, and every
FGN and LARCH short memory process. Finally, we realise a real application in the time series of hydroelectric dams in Uruguay. 
A summary of the main conclusions that can be drawn from the simulation study follows:

\begin{itemize}

\item The FOU test  has the best performance under the null hypothesis. There is more than one
family of short memory processes that includes several examples where the other tests drive to a wrong decision, while the 
FOU test does not make a mistake in its decision.
\item While the other tests can be wrong with a probability higher than desired  under both the null hypothesis and  the alternative hypothesis, the test proposed 
in this work is wrong with an excessively low probability under the null hypothesis. Therefore, when the null hypothesis is rejected in the FOU test, we can have greater
confidence that the observed time series really does have a long memory.
\item Under the  ARFIMA$(1,d,1)$ for $d>0$  where the parameters AR and MA are not very small, the FOU test has the best performance among all
of the tests considered, in terms of getting the best power.
\item In some cases, the FOU test in some cases is able to detect long memory processes when the other tests do not detect it.

\end{itemize}

\noindent \textbf{Acknowledgements}\\
We wish to thank Alejandro Cholaquidis for your help and support in the simulation study and Jos\'e Rafael
Le\'on for various rich conversations about this topic.

\end{document}